\documentclass[11pt]{amsart}
\usepackage[hmargin=2.5cm,vmargin=2.5cm]{geometry}
\usepackage{amsfonts}
\usepackage{amsthm}
\usepackage[english]{algorithm2e}
\usepackage{hyperref}
\usepackage{amsmath}
\usepackage{amssymb}
\usepackage[T1]{fontenc}
\newtheorem{Theorem}{Theorem}[section]
\newtheorem{Definition}[Theorem]{Definition}
\newtheorem{Remark}[Theorem]{Remark}
\newtheorem{Example}[Theorem]{Example}
\newtheorem{Lemma}[Theorem]{Lemma}
\newtheorem{Proposition}[Theorem]{Proposition}

\title{Enveloping Lie superalgebras and Killing-ricci forms of Bol superalgebras}
\author[Sylvain Attan]
       {Sylvain Attan
       \\\\
           D\'epartement de Math\'ematiques\\
       Universit\'e d'Abomey-calavi\\
       01 BP 4521, Cotonou 01, B\'enin\\
       syltane2010@yahoo.fr
       }
\begin{document}
\maketitle
\begin{abstract}
In this paper, enveloping Lie superalgebras of Bol superalgebras are introduced. The notion of Killing-Ricci forms and invariant forms of these superalgebras 
are investigated as generalization of the ones of Bol algebras. 
\end{abstract}
{\bf Mathematics Subject Classification:} 16W55, 17A30, 17A42, 17A70 

{\bf Keywords:} Bol superalgebras, enveloping Lie superalgebras, Killing-Ricci forms, invariant forms.
\section{Introduction}
Bol algebras were introduced in \cite{Mikheev1} in the context of a study of smooth Bol loops defined by the identity
\begin{eqnarray}
(xy\cdot z)y=x(yx\cdot y).\nonumber
\end{eqnarray}
The algebras play the same role with respect to Bol loops as Lie algebras do with respect to Lie groups  or Malt'sev algebras, to Moufang loops \cite{Sabinin1}.

A vector space V equipped with a trilinear operation $[,,]$ is called a Lie triple system  if
$$[a, a, b] = 0,$$
$$[a, b, c] + [b, c, a] + [c, a, b] = 0,$$
$$[x, y, [a, b, c] = [[x, y, a], b, c] +[ a, [x, y, b], c] +[ a, b, [x, y, c]]$$
for all $x, y, a, b, c \in V.$ A (left) Bol algebra $(V , [ , , ], [ , ])$ is a Lie triple system $(V , [ , , ])$
with an additional bilinear skew-symmetric operation $\cdot$ satisfying
$$[a, b, c\cdot d] = [a, b, c]\cdot d + c\cdot[a, b, d] + [c, d, a\cdot b] + 
(a\cdot b)\cdot(c\cdot d).$$
A related notion is that of a Lie triple algebra, introduced under the name generalized Lie triple systems, by Yamaguti \cite{Yamaguti1} and called later Lie Yamaguti algebra \cite{Kinyon1}.

From the standard enveloping Lie algebra of a given Bol algebra, the notion of Killing-Ricci form and invariant form for a Bol algebra are introduced and studied in \cite{Kuzmin}. 

A $\mathbb{Z}_2$ graded generalization of Lie algebras, called Lie superalgebras 
is considered in \cite{Kac, Scheunert} while a $\mathbb{Z}_2$ graded generalization of Lie Yamaguti algebras called Lie Yamaguti superalgebras were first considered 
 in \cite{Ouedraogo} and generalize Lie supertriple systems \cite{Tilgner}. The reader may refer to \cite{Okubo} for applications of Lie supertriple systems in physics.
 As Lie Yamaguti superalgebras, Bol superalgebras first introduced in \cite{rukavi} may also be viewed as a generalization of Lie supertriple systems. For  relations between Malcev superalgebras and Bol superalgebras, one may refer to \cite{attan}.
 
 As a part of the general theory of superalgebras, the notion of Killing form of Lie algebras is extended to the one of Lie triple systems \cite{Ravisankar}, Lie superalgebras \cite{Scheunert}, Lie supertriple systems \cite{Kamiya}, \cite{Peipei}, and next Lie Yamaguti superalgebras \cite{Zoungrana}.
  
  In this paper, we define enveloping Lie superalgebras and the Killing-Ricci form of Bol superalgebras, study this Killing-Ricci forms which could be seen as a generalization of the one of Bol algebras \cite{Kuzmin}  and the Killing form of Lie supertriple systems \cite{Kamiya}, \cite{Peipei}.
   
  The rest of this paper is organized as follows.\\
In section 2, we first recall some basics on Lie and Malcev superalgebras as well as Lie supertriple systems and Bol superalgebras. In section 3, we define the 
notion of pseudo superderivations (Definition \ref{Defpseudoder} ), study their properties (Lemma \ref{lpseudo}) and introduce the notion of enveloping Lie superalgebras of Bol superalgebras
(Definition \ref{DefEnveloping} ). In section 4, the Killing-Ricci form of Bol superalgebras is defined (Definition \ref{DefKR}) and some of its properties are investigated (Theorem \ref{Ricci}, Proposition \ref{consistence}, Lemma \ref{invar}). In the next section, the 
invariant form of Bol superalgebras is defined (Definition \ref{DefInvar}) and 
some results are obtained ( Lemma \ref{Iperp} , Theorem \ref{ThIperp}).

Throughout this paper, all vector superspaces and superalgebras are finite dimensional over a fixed 
ground field $\mathbb{K}$ of characteristic $0.$
\section{Some basics on superalgebras}
We recall here some useful definitions and examples of Lie supertriple systems as well as the ones of Bol superalgebras. These examples are obtained from the relation between Malcev superalgebras and Bol superalgebras which could be found in \cite{attan}.\\
Now let $M$ be a linear super-space over $\mathbb{K}$ that is a $\mathbb{Z}_2$-graded linear space with a direct sum $M=M_0\oplus M_1.$ 
The elements of $M_j,$ $j\in \mathbb{Z}_2,$ are said to be homogeneous of parity $j.$ The parity of a homogeneous element $x$ is denoted by $\bar{x}.$ 
For all $i,j\in \mathbb{Z}_2,$ $i+j$ will always means that this sum is calculated modulo $2.$ If $N=N_0\oplus N_1$ is another superspace, a linear map $f: M\rightarrow N$ is said to be of degree $r\in\mathbb{Z}_2$ if 
$f(M_i)\subseteq N_{i+r}$ for all $i\in\mathbb{Z}_2.$  If $f$ is of degree $r=0,$ that is $f(M_i)\subseteq N_i$ for all $i\in\mathbb{Z}_2,$ then $f$ is said to be an even linear map. An algebra $(A,[,])$ is called a superalgebra if the underlying vector space is $\mathbb{Z}_2$-graded i.e $A=A_0\oplus A_1$ and if furthermore $[A_i,A_j]\subset A_{i+j}.$
\begin{Definition}
 A Lie superalgebra is the superalgebra ($A=A_0\oplus A_1, [,]$)\\ satisfying
the super skew-symmetry and the super Jacobi identities that is
\begin{eqnarray}
[x,y]=-(-1)^{\bar{x}\bar{y}}[y,x]  
\end{eqnarray}
\begin{eqnarray}
 [[x,y],z]+(-1)^{\bar{x}(\bar{y}+\bar{z})}[[y,z],x]+(-1)^{\bar{z}(\bar{x}+\bar{y})}[[z,x],y]=0\label{superJacobi}
\end{eqnarray}
 for all $x,y,z\in \mathcal{H}(A).$ In term of the super Jacobian 
\begin{eqnarray}
 SJ(x,y,z)=[[x,y],z]+(-1)^{\bar{x}(\bar{y}+\bar{z})}[[y,z],x]+(-1)^{\bar{z}(\bar{x}+\bar{y})}[[z,x],y]
\end{eqnarray}
the superJacobi identity is written $SJ(x,y,z)=0$
\end{Definition}
for all $x,y,z\in\mathcal{H}(A).$\\

Another class of superalgebras that is of interest in this paper is the one of Malcev superalgebras.
\begin{Definition} \cite{albu1}\cite{alb2} A superalgebra $(M = M_0\oplus M_1,[,])$ is
called a Malcev superalgebra if it satisfies the following super identities:
 \begin{eqnarray}
 && [x,y]=-(-1)^{\bar{x}\bar{y}}[y,x]   \mbox{  (super skew-symmetry),} \nonumber\\
  &&[[[x,y],z],t]-[x,[[y,z],t]]-(-1)^{\bar{y}(\bar{z}+\bar{t})}[[x,[z,t]],y]-(-1)^{\bar{t}(\bar{y}+\bar{z})}[[[x,t],y],z]\nonumber\\
&=&(-1)^{\bar{y}\bar{z}}[[[ x,z],[y,t]]   \mbox{  (super Malcev identity)} \label{supermalcev}
 \end{eqnarray}
for  $x,y,z,t\in \mathcal{H}(M).$
\end{Definition}
\begin{Definition}(1) A supertriple system is a couple $(S=S_0\oplus S_1, [,,])$  consisting of a $\mathbb{Z}_2$-graded $\mathbb{K}$-vector space $S=S_0\oplus 
S_1$ and a $\mathbb{K}$-trilinear map $[,,]$ satisfying  
 $[S_i,S_j,S_k]\subset S_{i+j+k}$ for all $i,j,k$ in $\mathbb{Z}_2$  such that for all $x,y, z,\in \mathcal{H}(S)$ the following hold
\begin{eqnarray}
 [x,y,z]=-(-1)^{\bar{x}\bar{y}}[y,x,z] \ (\mbox{left super skew-symmetry})\nonumber
\end{eqnarray}
\begin{eqnarray}
 &&[x,y,z]+(-1)^{\bar{x}(\bar{y}+\bar{z})}[y,z,x]+(-1)^{\bar{z}(\bar{x}+\bar{y})}[z,x,y]=0\nonumber\\
 && \mbox{ (super ternary Jacobi identity)}\nonumber
\end{eqnarray}
(2) A  Lie supertriple system \cite{Tilgner},\cite{plz} is a supertriple system $(S=S_0\oplus S_1, [,,])$ such that the super ternary Nambu identity 
\begin{eqnarray}
 [x,y,[u,v,w]]&=&[[x,y,u],v,w]+(-1)^{\bar{u}(\bar{x}+\bar{y})}[u,[x,y,v],w]\nonumber\\
       &&+(-1)^{(\bar{x}+\bar{y})(\bar{u}+\bar{v})}[u,v,[x,y,w]]\label{supernambu}
\end{eqnarray}
 holds  for all $x,y, u, v, w \in \mathcal{H}(S).$
\end{Definition}
\begin{Example}\cite{attan, Tilgner, plz} Let $(L= L_0\oplus L_1, [,])$ be a Lie superalgebra. Then $(L,[,,])$ is a Lie supertriple system where for all $x,y, z\in \mathcal{H}(L),$ $[x,y,z]:= [[x,y],z].$
\end{Example}
\begin{Definition}\cite{attan} \cite{rukavi}
A Bol superalgebra is a triple $(B=B_0\oplus B_1,\cdot, [,,])$ consisting of a sperspace $B=B_0\oplus B_1,$ a linear map $\cdot: B^{\otimes 2}\rightarrow B$ 
satisfying $B_i\cdot B_j\subseteq B_{i+j}$ 
and a trilinear map $[,,]: B^{\otimes 3}\rightarrow B $ satisfying  $[B_i,B_j,B_k]\subseteq B_{i+j+k},$ $i, j, k \in\mathbb{Z}_2,$  such that: 
\\
$(BS_1)$ $x\cdot y= -(-1)^{\bar{x}\bar{y}} y\cdot x,$\\
$(BS_2)$ $[x,y,z]=-(-1)^{\bar{x}\bar{y}}[y,x,z],$\\
 $(BS_3)$ $[x,y,z]+(-1)^{\bar{x}(\bar{y}+\bar{z})}[y,z,x]+(-1)^{\bar{z}(\bar{x}+\bar{y})}[z,x,y]=0,$\\
$(BS_4)$ $[x,y,[u,v,w]]=[[x,y,u],v,w]+(-1)^{\bar{u}(\bar{x}+\bar{y})}[u,[x,y,v],w]+(-1)^{(\bar{x}+\bar{y})(\bar{u}+\bar{v})}[u,v,[x,y,w]],$\\
$(BS_5)$ $[x,y,u\cdot v]= (-1)^{\bar{u}(\bar{x}+\bar{y})}u\cdot[x,y,v]+[x,y,u]\cdot v+(-1)^{(\bar{x}+\bar{y})(\bar{u}+\bar{v})}[u,v ,x\cdot y]+ (x\cdot y)\cdot (u\cdot v)$\\
 holds  for all $x,y, u, v, w \in \mathcal{H}(S).$
\end{Definition}
\begin{Remark} A Bol superalgebra is a Lie triple supersystem $(B=B_0\oplus B_1, [,,])$ with a super anticommutative binary map $[,]$ such that
\begin{eqnarray}
 [x,y,u\cdot v]&= &(-1)^{\bar{u}(\bar{x}+\bar{y})}u\cdot[x,y,v]+[x,y,u]\cdot v+(-1)^{(\bar{x}+\bar{y})(\bar{u}+\bar{v})}[u,v ,x\cdot y]\nonumber\\
    &&+(x\cdot y)\cdot (u\cdot v)\label{ipseudoderivation}
\end{eqnarray}
holds forall $x, y, u, v\in \mathcal{H}(B).$
\end{Remark}
In \cite{attan} we proved that any Malcev superalgebra ($M = M_0\oplus M_1, [,]$)  equipped with a trilinear operation $\{,,\}$ where
\begin{eqnarray}
\{x,y,z\}&:=&\frac{1}{3}(2[[x,y],z]-(-1)^{\bar{x}(\bar{y}+\bar{z})}[[y,z],x]-(-1)^{\bar{z}(\bar{x}+\bar{y})}[[z,x],y])\nonumber
\end{eqnarray}
 for all $x, y, z\in \mathcal{H}(M),$ becomes a Bol superalgebra 
 $(M= M_0\oplus M_1, [,], \{,,\} ).$ This allows us to get the following
  example of Bol superalgebra from an example of a Malcev superalgebra \cite{albu3}.
\begin{Example}
  Let $(L^2(2,2), [,])$ be a non-Lie Malcev superalgebra \cite{albu3} defined with respect to a basis $(e_1,e_2,e_3,e_4)$ where $L^2(2,2)_{\bar{0}}=span(e_1,e_2),$ $L^2(2,2)_{\bar{1}}=span(e_3,e_4)$ with the nonzero products given by $[e_1,e_2]=e_2,$ $[e_1,e_3]=e_3,\ [e_1,e_4]=-e_4,\ [e_2,e_3]=-e_4$. Then $(L^2(3,1), [,],\{,,\})$ is a four-dimensional Bol superalgebra with
 $L^2(2,2)_{\bar{0}}=span(e_1,e_2),$ $L^2(2,2)_{\bar{1}}=span(e_3,e_4),$ where the nonzero products are given by
 $[e_1,e_2]=e_2,$ $[e_1,e_3]=e_3,\ [e_1,e_4]=-e_4,\ [e_2,e_3]=-e_4$ and\\
 $\{e_1,e_2,e_1\}=-e_2,$ $\{e_1,e_3,e_1\}=-e_3,$ $\{e_1,e_4,e_1\}=-e_4.$
 \hfill $\square$
\end{Example}
From \cite{attan} we  also get the following example.
 \begin{Example}\cite{attan} There is a four-dimensional Bol superalgebra
 $(L^2(3,1), [,],\{,,\})$  with\\
 $L^2(3,1)_{\bar{0}}=span(e_1,e_2,e_3),$ $L^2(3,1)_{\bar{1}}=span(e_4),$ where the nonzero products are given by\\
 $[e_1,e_3]=e_1,$ $[e_2,e_3]=e_1+e_2,$ $[e_3,e_4]=e_4, \ [e_4,e_4]=e_1$ and
 $\{e_1,e_3,e_3\}=e_1,$ $\{e_2,e_3,e_3\}=2e_1+e_2,$ $\{e_3,e_4,e_3\}=-e_4.$
\end{Example}
\begin{Definition}
Let $B=B_0\oplus B_1$ be a Bol superalgebra. A graded subspace $H=H_0\oplus H_1$ of 
$B$ is a subsuperalgebra of $B$ if $H_i\cdot H_j\subseteq H_{i+j}$ 
and $[H_i,H_j,H_k]\subseteq H_{i+j+k},$ for all $i, j, k \in\mathbb{Z}_2.$ 
\hfill $\square$
\end{Definition}
\begin{Definition}
A subsuperalgebra $H=H_0\oplus H_1$ of a Bol superalgebra $B$ is an 
invariant subsuperalgebra (resp., an ideal) of $B$ if $[B,B,H]\subseteq H$(resp.,
$BH\subseteq H$ and $[B,B,H]\subseteq H$). \hfill $\square$
\end{Definition}
If $H$ is an ideal of $B,$ it is an invariant subsuperalgebra of $B.$ Obvious 
the center of a Bol superalgebra $B$ defined by 
$Z(B)=\{x\in B, xy=0\ and\ [x,y,z]=0, \forall y, z\in B\}$ is an ideal of $B.$
\begin{Definition}
Let $A=A_0\oplus A_1$ and $B=B_0\oplus B_1$ be two Bol superalgebras. An even 
linear map $f: A \rightarrow B$ is called a morphism of Bol superalgebras if 
 $f(xy)=f(x)f(y)$ and $f([x,y,z])=[f(x),f(y),f(z)]$ for all $x,y,z\in B_0\cup B_1.$ \hfill $\square$
\end{Definition}
Recall \cite{Kac} that if $V=V_0\oplus V_1$ is a vector
superspace then, the set of the linear mappings of V into
itself which are homogeneous of degree $r$ is denoted by $End_r(V)=\{f \in End(V), f(V_i)\subseteq V_{r+i}\},$ we obtain an associative superalgebra $End(V)=End_0(V)\oplus End_1(V).$ The bracket
 $[f,g] =fg-(-1)^{\bar{f}\bar{g}}gf$ makes $End(V$) into a Lie superalgebra which
we denote by $l(V)$ or $l(m,n)$ where $m=dim V_0$ and 
$n=dim V_1.$ Let $e_1, \cdots, e_m, e_{m+1}, \cdots, e_{m+n}$ be a basis of $V.$ 
In this basis the matrix  of $f\in l(m,n)$ is expressed as 
$\left(\begin{array}{rl}
\alpha\ \beta\\
\gamma\ \delta
\end{array}
\right),$ 
$\alpha$ being an $(m\times m)-$,  $\delta$ an $(n\times n)-,$ $\beta$ an $(m\times n)-,$ and $\gamma$ an $(n\times m)-$ matrix. The matrices of even elements have the form $\left(\begin{array}{rl}
\alpha\ 0\\
0\ \delta
\end{array}
\right)$
and those of odd ones $\left(\begin{array}{rl}
0\ \beta\\
\gamma\ 0
\end{array}
\right)$. For $f=\left(\begin{array}{rl}
\alpha\ \beta\\
\gamma\ \delta
\end{array}
\right)$, the supertrace of
$M$ is defined by $str(M) = tr\alpha-tr\delta$ and does not depend on the
choice of a homogeneous basis. We have $str([f, g])=0$ that is
$str(fg)=(-1)^{\bar{f}\bar{g}}str(gf)$ and $str(hfh^{-1})=str(f).$
\section{Enveloping Lie superalgebras of a Bol superalgebra}
As derivations for algebras, superderivations of different superalgebras are an important subject of study in superalgebra and diverse area. They appear in many fields of mathematics and physics. In particular, they allow the construction of new superalgebras structures. In the case of Bol superalgebras, instead of superderivations, we have the notion of pseudo superderivations. They generalize  pseudo derivations for Bol algebras \cite{Sabinin1}  and superderivations for Lie supertriple systems \cite{plz} and allow the construction of enveloping Lie 
superalgebras of Bol superalgebras.
\begin{Definition}\label{Defpseudoder}
Let $B=B_0\oplus B_1$ be a Bol superalgebra. A linear map $P\in End_r(B)$ is called a pseudo superderivation of companion $a\in B_r ($ $r\in\mathbb{Z}_2$) if, for any $x,y,z\in B_0\cup B_1,$
 \begin{eqnarray}
 P([x,y,z]&=& [P(x),y,z]+(-1)^{r\bar{x}}[x,P(y),z]+(-1)^{r(\bar{x}+\bar{y})}[x,y,P(z)]\label{pseudo1}\\
 P(xy)&=&(-1)^{r\bar{x}}xP(y)+P(x)y+(-1)^{r(\bar{x}+\bar{y})}[x,y,a]+ a\cdot xy\label{pseudo2}  \hspace*{4cm}\mbox{   \hfill $\square$}
 \end{eqnarray}
\end{Definition}
\begin{Remark}
Note that a pseudo superderivation $P$ can have more than one companion. Let denote the set of all companions of a pseudo superderivation $P$ by $Com(P)$. 
\end{Remark}
Let $pS_r(B)$ be the set of all pseudo superderivations of degree $r$ and 
$pS(B)=pS_0(B)\oplus pS_1(B).$ Further let $PS_r(B)=\{ (P, a), P\in pS_r(B), a\in Com(P)\}$ and $PS(B)=PS_0(B)\oplus PS_1(B).$\\

For any pseudo superderivations $P,Q\in pS_0(B)\cup pS_1(B)$ of a Bol superalgebra $B,$ consider the supercommutator given by $[P,Q]:=PQ-(-1)^{\bar{P}\bar{Q}}QP$. Then we have:
\begin{Lemma}\label{lpseudo}
Let $(P,a)\in PS_r(B),$ $(Q,b)\in PS_r(B),$ $(R,c)\in PS_s(B)$ and $\lambda\in \mathbb{K}.$ Then
\begin{enumerate}
\item $(P,a)+\lambda(Q,b):=(P+\lambda Q,a+\lambda b)\in PS_r(B).$
\item $[(P,a),(R,c)]:=([P,R], P(c)-(-1)^{rs}R(a)-ac)\in PS_{r+s}(B).$
\end{enumerate}
\end{Lemma}
{\bf Proof.} Let $(P,a)\in PS_r(B),$ $(Q,b)\in PS_r(B),$ $(R,c)\in PS_s(B)$ and $\lambda\in \mathbb{K}.$ \\
 The first statement  is straightforward computation.For the second statement, pick  \\ $x,y,z,u,v\in B_0\cup B_1.$ Then using repeatedly (\ref{pseudo1}) for 
 $P$ and $R,$ we prove (\ref{pseudo1})  for $[P,R]$ as follows
\begin{eqnarray}
[P,R]([x,y,w])&=&PR([x,y,w])-(-1)^{rs}RP ([x,y,w])\nonumber\\
&=&P([R(x),y,z]+(-1)^{s\bar{x}}[x,R(y),z]+(-1)^{s(\bar{x}+\bar{y})}[x,y,R(z)])\nonumber\\
&&-(-1)^{rs}R([P(x),y,z]+(-1)^{r\bar{x}}[x,P(y),z]+(-1)^{r(\bar{x}+\bar{y})}[x,y,P(z)])\nonumber\\
&=&[PR(x),y,z]+(-1)^{r(\bar{x}+s)}[R(x),P(y),z]+(-1)^{r(s+\bar{x}+\bar{y})}[R(x),y,P(z)]\nonumber\\
&&+(-1)^{s\bar{x}}[P(x),R(y),z]+(-1)^{(r+s)\bar{x}}[x,PR(y),z]\nonumber\\
&&+(-1)^{s\bar{x}+r(\bar{x}+\bar{y}+s)}[x,R(y),P(z)]
+(-1)^{s(\bar{x}+\bar{y})}[P(x),y,R(z)])\nonumber\\
&&+(-1)^{s(\bar{x}+\bar{y})+r\bar{x}}[x,P(y),R(z)])+(-1)^{(r+s)(\bar{x}+\bar{y})}[x,y,PR(z)])\nonumber\\
&&-(-1)^{rs}\{[RP(x),y,z]+(-1)^{s(\bar{x}+r)}[P(x),R(y),z]+(-1)^{s(r+\bar{x}+\bar{y})}[P(x),y,R(z)]\nonumber\\
&&+(-1)^{r\bar{x}}[R(x),P(y),z]+(-1)^{(r+s)\bar{x}}[x,RP(y),z]\nonumber\\
&&+(-1)^{r\bar{x}+s(\bar{x}+\bar{y}+r)}[x,P(y),R(z)]
+(-1)^{r(\bar{x}+\bar{y})}[R(x),y,P(z)])\nonumber\\
&&+(-1)^{r(\bar{x}+\bar{y})+s\bar{x}}[x,R(y),P(z)])+(-1)^{(r+s)(\bar{x}+\bar{y})}[x,y,RP(z)])\}\nonumber\\
&=&[[P,R](x),y,z]+(-1)^{(r+s)\bar{x}}[x,[P,R](y),z]+(-1)^{(r+s)(\bar{x}+\bar{y})}[x,y,[P,R](z)]\nonumber
\end{eqnarray}
Thus we get (\ref{pseudo1}) for $[P,R].$ To prove, (\ref{pseudo2}) for $[P,R],$
 note that if we use repeatedly (\ref{pseudo2}) for $P$ we get 
\begin{eqnarray}
P(c\cdot xy)&=&(-1)^{rs}c\cdot P(xy)+P(c)\cdot xy+(-1)^{r(\bar{x}+\bar{y}+s}[c,xy,a]+a\cdot(c\cdot xy)  \nonumber\\
 &=&+ (-1)^{r(\bar{x}+s)}c\cdot xP(y)+(-1)^{rs}c\cdot P(x)y+ (-1)^{r(\bar{x}+\bar{y}+s)}c\cdot[x,y,a]\nonumber\\
 &&+(-1)^{rs}c\cdot(a\cdot xy) +P(c)\cdot xy+(-1)^{r(\bar{x}+\bar{y}+s}[c,xy,a]+a\cdot(c\cdot xy)\label{pseudo3}
\end{eqnarray} 
  and therefore by (\ref{pseudo1}) and (\ref{pseudo3}) for $P$ and $R$ we get:
 \begin{eqnarray}
 [P,R](xy)&=&PR(xy)-(-1)^{rs}RP (xy)\nonumber\\
 &=&P((-1)^{s\bar{x}}xR(y)+R(x)y+(-1)^{s(\bar{x}+\bar{y})}[x,y,c]+ c\cdot xy)\nonumber\\
 &&-(-1)^{rs}\{R((-1)^{r\bar{x}}xP(y)+P(x)y+(-1)^{r(\bar{x}+\bar{y})}[x,y,a]+ a\cdot xy)\}\nonumber\\
 &=&(-1)^{(r+s)\bar{x}}xPR(y)+(-1)^{s\bar{x}}P(x)R(y)+(-1)^{r(\bar{x}+\bar{y}+s)+
 s\bar{x}}[x,R(y),a]+ (-1)^{s\bar{x}}a\cdot xR(y)      \nonumber\\
 &&+(-1)^{r(\bar{x}+s)}R(x)Py+PR(x)P(y)+(-1)^{r(\bar{x}+\bar{y}+s)}[R(x),y,a]+a\cdot R(x)y                    \nonumber\\
 &&+(-1)^{s(\bar{x}+\bar{y})}[P(x),y,c] +(-1)^{s(\bar{x}+\bar{y})+r\bar{x}}[x,P(y),c]+(-1)^{(r+s)(\bar{x}+\bar{y})}[x,y,P(c)]\nonumber\\
 && + (-1)^{r(\bar{x}+s)}c\cdot xP(y)+(-1)^{rs}c\cdot P(x)y+ (-1)^{r(\bar{x}+\bar{y}+s)}c\cdot[x,y,a]\nonumber\\
 &&+(-1)^{rs}c\cdot(a\cdot xy) +P(c)\cdot xy+(-1)^{r(\bar{x}+\bar{y}+s)}[c,xy,a]+a\cdot(c\cdot xy)\nonumber\\
 &&-(-1)^{rs}\{ (-1)^{(r+s)\bar{x}}xRP(y)+(-1)^{r\bar{x}}R(x)P(y)+(-1)^{s(\bar{x}+\bar{y}+r)+
 r\bar{x}}[x,P(y),c]      \nonumber\\
 &&+ (-1)^{r\bar{x}}c\cdot xP(y)+(-1)^{s(\bar{x}+r)}P(x)Ry+RP(x)(y)+(-1)^{s(\bar{x}+\bar{y}+r)}[P(x),y,c] \nonumber\\
 &&+c\cdot P(x)y+(-1)^{r(\bar{x}+\bar{y})}[R(x),y,a] +(-1)^{r(\bar{x}+\bar{y})+s\bar{x}}[x,R(y),a]\nonumber\\
 &&+(-1)^{(r+s)(\bar{x}+\bar{y})}[x,y,R(a)] + (-1)^{s(\bar{x}+r)}a\cdot xR(y)+(-1)^{rs}a\cdot R(x)y\nonumber\\
 &&+(-1)^{s(\bar{x}+\bar{y}+r)}a\cdot[x,y,c]+(-1)^{rs}a\cdot(c\cdot xy) +R(a)\cdot xy+(-1)^{s(\bar{x}+\bar{y}+r)}[a,xy,c]\nonumber\\
 &&+c\cdot(a\cdot xy)\}\nonumber
 \end{eqnarray}
\begin{eqnarray}
 &=&(-1)^{(r+s)\bar{x}}x\cdot[P,R](y)+[P,R](x)\cdot y+
 (-1)^{(r+s)(\bar{x}+\bar{y})}[x,y,P(c)-(-1)^{rs}R(a)]\nonumber\\
 &&+(P(c)-(-1)^{rs}R(a))\cdot xy+(-1)^{r(\bar{x}+\bar{y}+s)}[c,xy,a]-
 (-1)^{s(\bar{x}+\bar{y})}[a,xy,c]\nonumber\\
 &&+(-1)^{r(\bar{x}+\bar{y}+s)}c[x,y,a]-(-1)^{s(\bar{x}+\bar{y})}a[x,y,c]\nonumber\\
 &=&(-1)^{(r+s)\bar{x}}x\cdot[P,R](y)+[P,R](x)\cdot y+
 (-1)^{(r+s)(\bar{x}+\bar{y})}[x,y,P(c)-(-1)^{rs}R(a)]\nonumber\\
 &&+(P(c)-(-1)^{rs}R(a))\cdot xy-(-1)^{(r+s)(\bar{x}+\bar{y})+rs}[xy,c,a]-
 (-1)^{s(\bar{x}+\bar{y})}[a,xy,c]\nonumber\\
 &&-(-1)^{(r+s)(\bar{x}+\bar{y})}[x,y,a]c-(-1)^{s(\bar{x}+\bar{y})}a[x,y,c] 
 \mbox{ (by $(BS_1)$ and $(BS_2)$)}\nonumber\\
 &=&(-1)^{(r+s)\bar{x}}x\cdot[P,R](y)+[P,R](x)\cdot y+
 (-1)^{(r+s)(\bar{x}+\bar{y})}[x,y,P(c)-(-1)^{rs}R(a)-ac]\nonumber\\
 &&+(P(c)-(-1)^{rs}R(a)-ac)\cdot xy
 \mbox{ (by $(BS_1)-(BS_3)$ and $(BS_5)$)}\nonumber
 \end{eqnarray}
 then we get (\ref{pseudo2}). That ends the proof. \hfill $\square$\\
 
 Now, we can prove the following result.
\begin{Proposition}
Let $B=B_0\oplus B_1$ be a Bol superalgebra. Then  $(PS(B),[,])$ is a Lie superalgebra. 
\end{Proposition} 
{\bf Proof.} The proof follows by Lemma \ref{lpseudo} and the fact that 
$pS(B)$ is a Lie superalgebra since it is straightforward to check that $pS(B)$ is a subsuperalgebra of the Lie superalgebra $End(B).$ \hfill $\square$
\begin{Definition}
$PS(B)$ is called the enlarged Lie superalgebra of pseudo superderivations of a  Bol superalgebra $B.$ \hfill $\square$
\end{Definition}
Let $B=B_0\oplus B_1$ be a Bol superalgebra. For any $x,y\in B_0\cup B_1,$ denote by $D_{x,y},$ the endomorphism of $B$ defined by $D_{x,y}(z):=[x,y,z]$ for all
$z\in B.$ We have, for any $x,y\in B_0\cup B_1,$ $r\in\mathbb{Z}_2,$ $D_{x,y}(B_r)\in B_{r+\bar{x}+\bar{y}},$ that is, $D_{x,y}$ is a linear map of degree 
$\bar{x}+\bar{y}.$ Moreover, it comes from $(BS_4)$ and $(BS_5)$ that
\begin{eqnarray}
D_{x,y}([u,v,w])&=&[D_{x,y}(u),v,w]+(-1)^{\bar{u}(\bar{x}+\bar{y})}[u,D_{x,y}(v),w]\nonumber\\
 &&+(-1)^{(\bar{x}+\bar{y})(\bar{u}+\bar{v})}[u,v,D_{x,y}(w)]\\
D_{x,y}(u\cdot v)&=& (-1)^{\bar{u}(\bar{x}+\bar{y})}u\cdot D_{x,y}(v)+D_{x,y}(u)\cdot v+(-1)^{(\bar{x}+\bar{y})(\bar{u}+\bar{v})}[u,v ,x\cdot y]\nonumber\\
 &&+ (x\cdot y)\cdot (u\cdot v).
\end{eqnarray}
for any $x,y,u,v$ in $B_0\cup B_1.$ It follows that $D_{x,y}$ is pseudo superderivation of degree $\bar{x}+\bar{y}$ and companion $xy,$ called inner 
pseudo superderivation of $B.$ \\

Now one can reformulate the definition of a  Bol superalgebra in the following
manner.
\begin{Definition}
A vector superspace $B=B_0\oplus B_1$ equipped with a bilinear operation 
$(x,y)\longmapsto x\cdot y$ satisfying $B_i\cdot B_j\subseteq B_{i+j}$ 
and a trilinear operation $(x,y,z)\longmapsto [x,y,z]$ satisfying  $[B_i,B_j,B_k]\subseteq B_{i+j+k};$ $i, j, k \in\mathbb{Z}_2,$  is called a  Bol superalgebra if
\\
$(BS_1)$ $x\cdot y= -(-1)^{\bar{x}\bar{y}} y\cdot x,$\\
$(BS_2)$ $[x,y,z]=-(-1)^{\bar{x}\bar{y}}[y,x,z],$\\
 $(BS_3)$ $[x,y,z]+(-1)^{\bar{x}(\bar{y}+\bar{z})}[y,z,x]+(-1)^{\bar{z}(\bar{x}+\bar{y})}[z,x,y]=0,$
and the endomorphism\\ $D_{x,y}: z\longmapsto [x,y,z]$ is its pseudo
superderivation with a companion $xy$ for all $x,y\in B_0\cup B_1$ and 
$z\in B_0\cup B_1.$ \hfill $\square$
\end{Definition}
Let $ipS_r(B,B)$ be the vector space spanned by all inner pseudo superderivations 
$D_{x,y}$ $(x,y\in B_0\cup B_1\ and\ \bar{x}+\bar{y}=r\in\mathbb{Z}_2).$\\
We can define naturally a $\mathbb{Z}_2-$gradation by setting 
$ipS(B,B)=ipS_0(B,B)\oplus ipS_1(B,B).$\\
Evidently, $ipS(B,B)$ is a subsuperalgebra of the Lie superalgebra $pS(B).$

Accordingly, $IPS(B,B)$ can be introduced as the set of all pairs $(P,c),$
where $P\in ipS(B,B),$ $c\in Com(P).$
Evidently, $IPS(B,B)$ is a subsuperalgebra of the Lie superalgebra $PS(B).$
\begin{Definition}
$IPS(B,B)$ is called the enlarged Lie superalgebra of inner pseudo\\ superderivations. \hfill $\square$
\end{Definition}
Let $B=B_0\oplus B_1$ be a Bol superalgebra, $H=H_0\oplus H_1$ be a subsuperalgebra of $PS(B)$ such that $IPS(B,B)\subseteq H.$ For $i\in\mathbb{Z}_2,$ let $\mathcal{L}_i(B)=B_i\oplus H_i,$ and define a new superbracket operation in $\mathcal{L}^H(B)=\mathcal{L}_0(B)\oplus \mathcal{L}_1(B)=B\oplus H$ as follows:
for any $x,y\in B_0\cup B_1,$ $(P,a),(Q,b)\in H_0\cup H_1,$
\begin{eqnarray}
 \left[ x,y\right] &:=& (D_{x,y}, xy)\label{li1}\\
  \left[ (P,a),x\right] &:=&-(-1)^{\bar{a}\bar{x}}\left[ x,(P,a)\right] := P(x)\label{li2}\\
  \left[ (P,a),(Q,b)\right] &:=&(\left[ P,Q\right], P(b)-(-1)^{\bar{a}\bar{b}}Q(a)-ab)\label{li3}
\end{eqnarray}
Then we have:
\begin{Proposition}
$(\mathcal{L}^H(B), [,])$ is a Lie superalgebra. 
\end{Proposition}
{\bf Proof.} It is clear that the operation $[,]$ is supersymmetric.
For the Jacobi superidentity, there is many cases to distinguish. First 
for all $x,y, u$ in $B_0\cup B_1,$ $(P,a),(Q,b)$ in $H_0\cup H_1,$ we get
\begin{eqnarray}
&&([[x,y],(P,a)]+(-1)^{\bar{x}(\bar{y}+\bar{a})}[[y,(P,a)],x]+(-1)^{\bar{a}(\bar{x}+\bar{y})}[[(P,a),x],y])(u)\nonumber\\
&=&[x,y,P(u)]-(-1)^{\bar{a}(\bar{x}+\bar{y})}P([x,y,u])-(-1)^{\bar{x}(\bar{a}+\bar{y})+\bar{a}\bar{y}}[P(y),x,u]+(-1)^{\bar{a}(\bar{x}+\bar{y})}[P(x),y,u]\nonumber\\
&=&-(-1)^{\bar{a}(\bar{x}+\bar{y})}(-(-1)^{\bar{a}(\bar{x}+\bar{y})}[x,y,P(u)]+P([x,y,u])-(-1)^{\bar{a}\bar{x}}x,[P(y),u]-[P(x),y,u])\nonumber\\
&=& 0 \mbox{ (by  (\ref{pseudo1}) )}\nonumber
\end{eqnarray}
i.e., $[[x,y],(P,a)]+(-1)^{\bar{x}(\bar{y}+\bar{a})}[[y,(P,a)],x]+(-1)^{\bar{a}(\bar{x}+\bar{y})}[[(P,a),x],y]=0$
\begin{eqnarray}
&&[[P,a),x],(Q,b)]+(-1)^{\bar{a}(\bar{x}+\bar{b})}[[(x,(Q,b)],(P,a)]+
(-1)^{\bar{b}(\bar{a}+\bar{x})}[[(Q,b),(P,a)],x]\nonumber\\
&=&PQ(x)-(-1)^{\bar{a}\bar{b}}QP(x)+(-1)^{\bar{a}(\bar{x}+\bar{b})}[Q(x),(P,a]
-(-1)^{\bar{x}\bar{b}}[P(x),(Q,b)]\nonumber\\
&=& PQ(x)-(-1)^{\bar{a}\bar{b}}QP(x)-PQ(x)+(-1)^{\bar{a}\bar{b}}QP(x)=0\nonumber
\end{eqnarray}
Next the other cases when the three elements are in $B_0\cup B_1$ or $H_0\cup H_1$ follow from $(HBS_3)$ and the fact that $H$ is a Lie superalgebra. \hfill $\square$
\begin{Definition}\label{DefEnveloping}
An enveloping superalgebra of a Bol superalgebra $B$ is a Lie superalgebra 
  $\mathcal{L}^H(B)$ defined above.Taking $H=PS(B)$ we obtain the maximal enveloping superalgebra, taking $H=IPS(B,B)$ we obtain the minimal (standard) enveloping superalgebra. \hfill $\square$
\end{Definition}
The following result will be used in the last section.
\begin{Lemma}\label{Idnvelop}
Let $K$ be an ideal of a Bol superalgebra $B.$ Then $\mathcal{K}=K\oplus IPS(B,K)$ is an ideal of the standard enveloping superalgebra  $L=B\oplus IPS(B,B).$
\end{Lemma}
{\bf Proof.} It suffices to prove that $[\mathcal{K},L]\subseteq\mathcal{K}$ 
which is a straightforward computation. \hfill $\square$
\section{killing-Ricci forms of Bol superalgebras.}
The definition of the Killing-Ricci form given here for Bol superalgebras 
stems from \cite{Kuzmin} where the Killing-Ricci form for Bol algebras is 
defined following \cite{Kikkawa} as the restriction of the Killing form of the 
standard enveloping Lie algebra of the given Bol algebra to this later.
Let $B=B_0\oplus B_1$ be a an $n$-dimensional Bol superalgebra and 
$L(B)=(B_0\oplus IPS_0(B,B))\oplus (B_1\oplus IPS_1(B,B)):=B\oplus IPS(B,B)$ its standard  enveloping Lie superalgebra. Let $\alpha$ be the Killing form of $L(B)$ 
and $\beta,$ the restriction of $\alpha$ to $B\times B.$ From \cite{Kuzmin},
 we introduce the following definition:
\begin{Definition}\label{DefKR}
The form $\beta$ is called the Killing-Ricci form of the Bol superalgebra $B.$ 
\hfill $\square$
\end{Definition}
For any $x,y\in B_0\cup B_1,$ define the endomorphism $R_{x,y}$ of the vector 
superspace $B$  by $R_{x,y}(z)=(-1)^{\bar{z}(\bar{x}+\bar{y})}[z,x,y]=(-1)^{\bar{z}(\bar{x}+\bar{y})}D_{z,x}(y)$ for all $z\in B_0\cup B_1.$ It is clear that, 
$R_{x,y}$ is of degree $\bar{x}+\bar{y}.$\\

The next result gives an explicit expression of $\beta.$
\begin{Theorem}\label{Ricci}
For every $x,y\in B_0\cup B_1,$ 
\begin{eqnarray}
\beta(x,y)&=& str(R_{x,y}+(-1)^{\bar{x}\bar{y}}R_{y,x})\label{killing}
\end{eqnarray}
\end{Theorem}
{\bf Proof.} Let $\{e_i\},$ $\{f_i\},$ $\{u_i\},$ $\{v_i\}$ be bases for 
$B_0, B_1,$ $IPS_0(B,B),$ $IPS_1(B,B),$ respectively. It suffices to prove 
(\ref{killing}) for all elements $x,\ y$ in the basis i.e., 
$\beta(e_i,e_j)=\alpha(e_i,e_j),$ $\beta(e_i,f_j)=\alpha(e_i,f_j),$ $\beta(f_i,f_j)=\alpha(f_i,f_j)$.  For these bases, we express 
the operations of $B$ and $ISP(B,B),$  using the tensor notation
 (i.e., repeated indices imply summation), as follows:
\begin{eqnarray}
&& D_{e_i,e_j}=R_{ij}^m u_m,\ D_{e_i,f_j}=S_{ij}^mv_m,\ D_{f_i,f_j}=T_{ij}^mu_m, 
[u_m,e_i]=u_m(e_i)=C_{mi}^je_j,\nonumber\\
 &&[v_m,e_i]=v_m(e_i)=H_{mi}^jf_j,\  [u_m,f_i]=u_m(f_i)=K_{mi}^jf_j,\         
 [v_m,f_i]=v_m(f_i)= L_{mi}^je_j\nonumber
\end{eqnarray}
Since $R_{e_i,f_j}(B_0)\subseteq B_1,$ $R_{e_i,f_j}(B_1)\subseteq B_0,$ we have 
$str(R_{e_i,f_j}+R_{f_j,e_i})=0$ and then $\beta(e_i,f_j)=0=\alpha(e_i,f_j)$
(consistency property of $\alpha$). Hence, it remains to show that
 $\beta(e_i,e_j)=\alpha(e_i,e_j),$ $\beta(f_i,f_j)=\alpha(f_i,f_j).$ The
  identities (\ref{li1}) and (\ref{li2}) imply the following:
  \begin{eqnarray}
&& L_{e_i}L_{e_j}(e_k)=[e_i,[e_j,e_k]]=[e_i, D_{e_j,e_k}]= D_{e_k,e_j}(e_i)
   = R_{kj}^mC_{mi}^te_t\nonumber\\
&&  L_{e_i}L_{e_j}(f_k)=[e_i,[e_j,f_k]]=[e_i, D_{e_j,f_k}]= -D_{e_j,f_k}(e_i)
   = -S_{jk}^mH_{mi}^tf_t\nonumber\\  
&& L_{e_i}L_{e_j}(u_m)=[e_i,[e_j,u_m]]=-[e_i, u_m(e_j)]= -C_{mj}^k[e_i,e_k]
   =-C_{mj}^kR_{ik}^tu_t\nonumber\\
&& L_{e_i}L_{e_j}(v_m)=[e_i,[e_j,v_m]]=-[e_i, v_m(e_j)]= -H_{mj}^k[e_i,e_k]
   =-H_{mj}^kS_{ik}^tv_t.\nonumber
  \end{eqnarray}
  Thus from the relation $\beta(e_i,e_j)=\alpha(e_i,e_j)=str(L_{e_i}L_{e_j}),$
  we get
  \begin{eqnarray}
  \beta(e_i,e_j)=R_{kj}^mC_{mi}^k+S_{jk}^mH_{mi}^k-C_{mj}^kR_{ik}^m+H_{mj}^kS_{ik}^m\label{k1}
  \end{eqnarray}
  From the other hand, we have 
  \begin{eqnarray}
  && R_{e_i,e_j}(e_k)=D_{e_k,e_i}(e_j)=-D_{e_i,e_k}(e_j)=-R_{ik}^mu_m(e_j)
  =-R_{ik}^mC_{mj}^te_t\nonumber\\
  && R_{e_i,e_j}(f_k)=D_{f_k,e_i}(e_j)=-D_{e_i,f_k}(e_j)=-S_{ik}^mv_m(e_j)
  =-S_{ik}^mH_{mj}^tf_t\nonumber
  \end{eqnarray}
  By interchanging $i$ and $j,$ we have
  \begin{eqnarray}
  R_{e_j,e_i}(e_k)
  =-R_{jk}^mC_{mi}^te_t=R_{kj}^mC_{mi}^te_t \mbox{ and }
   R_{e_j,e_i}(f_k)=-R_{jk}^mH_{mi}^tf_t.\nonumber
  \end{eqnarray}
  Then we get 
  \begin{eqnarray}
 str(R_{e_i,e_j}+(-1)^{\bar{ei}\bar{e_j}} R_{e_j,e_i})&=&-R_{ik}^mC_{mj}^k+S_{ik}^mH_{mj}^k+R_{kj}^mC_{mi}^k+
  R_{jk}^mH_{mi}^k\label{k2}
  \end{eqnarray}
  from (\ref{k1}) and (\ref{k2}), we obtain
  \begin{eqnarray}
  \beta(e_i,e_j)&=&str(R_{e_i,e_j}+(-1)^{\bar{ei}\bar{e_j}} R_{e_j,e_i})\nonumber
  \end{eqnarray}
  Again, the identities (\ref{li1}) and (\ref{li2}) imply the following:
  \begin{eqnarray}
&& L_{f_i}L_{f_j}(e_k)=[f_i,[f_j,e_k]]=-[f_i, D_{e_k,f_j}]=-D_{e_k,f_j}(f_i)=-S_{kj}^mv_m(f_i)=-S_{kj}^mL_{mi}^te_t\nonumber\\
&& L_{f_i}L_{f_j}(f_k)=[f_i,[f_j,f_k]]=[f_i, D_{f_j,f_k}]=-D_{f_j,f_k}(f_i)=-T_{jk}^mu_m(f_i)=-T_{jk}^mK_{mi}^tf_t\nonumber\\
&& L_{f_i}L_{f_j}(u_k)=[f_i,[f_j,u_k]]=-[f_i, u_k(f_j)]=-K_{kj}^m[f_i,f_m]
=-K_{kj}^mD_{f_i,f_m}=-K_{kj}^mT_{im}^su_s\nonumber\\
&& L_{f_i}L_{f_j}(v_k)=[f_i,[f_j,v_k]]=[f_i, v_k(f_j)]=L_{kj}^m[f_i,e_m]
=-L_{kj}^mD_{e_m,f_i}=-L_{kj}^mS_{mi}^sv_s\nonumber
  \end{eqnarray}
 Again from the relation $\beta(f_i,f_j)=\alpha(f_i,f_j)=str(L_{f_i}L_{f_j}),$
 we get
 \begin{eqnarray}
 \beta(f_i,f_j)=-S_{kj}^mL_{mi}^k+T_{jk}^mK_{mi}^k-K_{kj}^mT_{im}^k+
 L_{kj}^mS_{mi}^k\label{k3}
 \end{eqnarray}
 From the other hand, we have 
  \begin{eqnarray}
   R_{f_i,f_j}(e_k)=D_{e_k,f_i}(f_j)=S_{ki}^mv_m(f_j)=S_{ki}^mL_{mj}^se_s
  \mbox{ and }
   R_{f_i,f_j}(f_k)=D_{f_k,f_i}(f_j)=T_{ki}^mu_m(f_j)=T_{ki}^mK_{mj}^sf_s
  \nonumber
  \end{eqnarray}
  By interchanging $i$ and $j,$ we get
  \begin{eqnarray}
  R_{f_j,f_i}(e_k)=S_{kj}^mL_{mi}^se_s \mbox{ and }
  R_{f_j,f_i}(f_k)=T_{kj}^mK_{mi}^sf_s=T_{jk}^mK_{mi}^sf_s.\nonumber
  \end{eqnarray}
  Next we get
  \begin{eqnarray}
  str(R_{f_i,f_j}+(-1)^{\bar{fi}\bar{f_j}}R_{f_j,f_i})=str(R_{f_i,f_j}-R_{f_j,f_i})&=&S_{ki}^mL_{mj}^k-T_{ki}^mK_{mj}^k
  -S_{kj}^mL_{mi}^k+T_{jk}^mK_{mi}^k\nonumber\\
  &=&S_{mi}^kL_{kj}^m-T_{mi}^kK_{kj}^m
  -S_{kj}^mL_{mi}^k+T_{jk}^mK_{mi}^k\label{k4}
  \end{eqnarray}
  and therefore from (\ref{k3}) and (\ref{k4}), we obtain
  \begin{equation}
  \beta(f_i,f_j)=str(R_{f_i,f_j}+(-1)^{\bar{fi}\bar{f_j}}R_{f_j,f_i})\nonumber
  \end{equation}
  whereby relation (\ref{killing}) is proved. \hfill $\square$
  \begin{Remark}
  Recall that if $(B, \cdot, [,,])$ is a (left) Bol algebra, then the Killing-Ricci form on $B$ is 
  defined as $\beta(x,y)=tr(R_{x,y}+R_{y,x})$ which we deduce from the one for (right) Bol algebra \cite{Kuzmin} where $R_{x,y}(z)=[z,x,y].$ So if a Bol superalgebra $B$ is reduced to a Bol algebra, $\beta$ as in Theorem \ref{Ricci} is the Killing-Ricci form of the Bol algebra $B.$
  \end{Remark}
\begin{Proposition} \label{consistence}
  Let $B=B_0\oplus B_1$ be a Bol superalgebra with a Killing-Ricci form $\beta.$ Then
  \begin{enumerate}
  \item $\beta(x,y)=(-1)^{\bar{x}\bar{y}}\beta(y,x)$ for all $x,y\in B_0\cup B_1$        
       (supersymmetry),
  \item $\beta(B_0,B_1)=0 $ (consistence),
  \item $\beta(A(x),A(y))=\beta(x,y),$ $A\in Aut(B),$
  \item $\beta([x,y,z],u)=-(-1)^{\bar{z}(\bar{x}+\bar{y})}\beta(z,[x,y,u]),$
   for all $x, y, z, u\in B_0\cup B_1.$  
  \end{enumerate}
\end{Proposition}
{\bf Proof.} We know \cite{Kac} that if $\alpha$ is a Killing form of any 
Lie superalgebra, then $\alpha$ satisfies the relation $(1)-(3)$ of the 
Proposition above. Then the relation $(1)-(3)$ follows from the fact that 
the Killing-Ricci form $\beta$ is the restriction to $B$ of the Killing form of 
the standard enveloping superalgebra $L(B)$ of $B.$  For the relation (4), pick
$x, y, z, u\in B_0\cup B_1$ et denote by $\alpha$ the Killing form of $L(B).$
Since $L(B)$ is a Lie superalgebra, then $\alpha$ satisfies
\begin{eqnarray}
\alpha([x,y],z)=-(-1)^{\bar{x}\bar{y}}\alpha(y,[x,z])\nonumber
\end{eqnarray}
then 
\begin{equation}
\beta([x,y,z],u)=\alpha([[x,y],z],u)=-(-1)^{\bar{z}(\bar{x}+\bar{y})}\alpha(z,[[x,y],u])=-(-1)^{\bar{z}(\bar{x}+\bar{y})}\beta(z,[x,y,u]).\nonumber
\end{equation}
Hence the proposition is proved. \hfill $\square$
\begin{Remark}
If we consider the Bol superalgebra $B$  as a Lie supertriple system with the ternary operation $[,,],$ that we denote by $B_S,$ then the relation $(4)$ in Proposition \label{consistence}, says that the Killing-Ricci form of $B$ is an 
invariance form of $B_S.$
\end{Remark}
The following result will be used below.
\begin{Lemma}\label{invar}
Let $K$ be a Killing-Ricci form of a Bol superalgebra $B.$ Then the following 
conditions are equivalent:
\begin{eqnarray}
 K([x,y,z],u)&=&-(-1)^{\bar{z}(\bar{x}+\bar{y})}K(z,[x,y,u])\nonumber\\
 K([x,y,z],u)&=&-(-1)^{\bar{y}(\bar{z}+\bar{u})}K(x,[z,u,y])\nonumber\\
 K(x,[y,z,u])&=&(-1)^{\bar{x}\bar{y}+\bar{z}\bar{u}}K(y,[x,u,z])\nonumber
\end{eqnarray}
for all $x, y, z, u\in B_0\cup B_1.$
\end{Lemma}
{\bf Proof.} We know that in any Lie supertriple system with an invariant form,  the relations (\ref{inva1})-(\ref{inva3}) are inquivalent \cite{Kamiya}. The proof then follows from the fact that the Killing-Ricci form of $B$ is an invariant form of $B_S.$ \hfill $\square$
\section{Invariant forms of Bol superalgebras}
In this section, we introduce the concept of invariant forms of 
Bol superalgebras as generalization of those of Bol algebras and Lie supertriple 
systems.
\begin{Definition}\label{DefInvar}
An invariant form $b$ of a Bol superalgebra $B=B_0\oplus B_1$ is a 
supersymmetric bilinear form on $B$ satisfying the identities
\begin{eqnarray}
b(xy,z)&=&-(-1)^{\bar{x}\bar{y}}b(y,xz)\label{inv1}\\
b([x,y,z],u)&=&-(-1)^{\bar{y}(\bar{z}+\bar{u})}b(x,[z,u,y])\label{inv2}
\end{eqnarray}
for all $x, y, z, u\in B_0\cup B_1.$ \hfill $\square$
\end{Definition} 
\begin{Remark}
If $B$ is reduced to a Lie supertriple system (resp. a (left) Bol algebra), then 
$b$ is reduced to an invariant form to a Lie supertriple system \cite{Kamiya, Peipei}, (resp. a (left) 
Bol algebra) which can be deduced from the one of a (right) Bol algebra \cite{Kuzmin} 
\end{Remark}
Let $b$ be an invariant form of a Bol superalgebra $B=B_0\oplus B_1.$ Then $b$ is 
an invariant form of $B_S$ and by Lemma \ref{invar}, $b$ satifies the following 
equivalent conditions:
\begin{eqnarray}
b([x,y,z],u)&=&-(-1)^{\bar{z}(\bar{x}+\bar{y})}b(z,[x,y,u])\label{inva1}\\
 b([x,y,z],u)&=&-(-1)^{\bar{y}(\bar{z}+\bar{u})}b(x,[z,u,y])\label{inva2}\\
 b(x,[y,z,u])&=&(-1)^{\bar{x}\bar{y}+\bar{z}\bar{u}}b(y,[x,u,z])\label{inva3}
\end{eqnarray}
for all $x, y, z, u\in B_0\cup B_1.$
\begin{Definition}
Let $b$ be an invariant form of a Bol superalgebra $B$ and $V$ be a subset of $B.$
The orthogonal $V^\perp$ of $V$ with respect to $b$ is defined by 
$V^\perp=\{x\in B, b(x,y)=0,\ \forall y\in V\}.$ The invariant form $b$ 
is nondegenerate if $B^\perp=\{0\}.$ \hfill $\square$
\end{Definition}
\begin{Lemma}\label{Iperp}
Let $b$ be an invariant form of a Bol superalgebra $B.$ Then,
\begin{enumerate}
\item $(B+[B,B,B])^\perp=Z(B)$ if $b$ is nondegenerate;
\item  If $I$ is an ideal of $B$ then $I^\perp$ is an ideal of $B.$ In 
particular,  $B^\perp$ is  an ideal of $B.$
\end{enumerate}
\end{Lemma}
{\bf Proof.} Pick an homogeneous element $u$ in $(B+[B,B,B])^\perp.$ Then, for any $x,y,z$ in $B_0\cup B_1,$ we get $b(u,xy)=0$ and $b(u,[x,y,z])=0.$ This implies by (\ref{inv1})  and (\ref{inv2}) that $(-1)^{\bar{u}\bar{x}}b(xu,y)=0$ and   $(-1)^{\bar{z}(\bar{x}+\bar{y})}b([u,z,x],y)=0.$ As $b$ is nondegenerate, we obtain $xu=0$ and $[u,z,x]=0$ for any $x, z$ in $B_0\cup B_1,$ this gives $u\in Z(B).$ 

Conversely, if $u\in Z(B),$ we have for all $x, y, z, v, w$ in $B,$ 
$b(u, xy+[z,v,w])=b(u,xy)+b(u,[z,v,w])=0$ and $u\in (B+[B,B,B])^\perp$ whence 
$(B+[B,B,B])^\perp=Z(B).$

Now, suppose that $I$ is an ideal of $B$ that is $BI\subseteq I,$ 
$[B,I,B]\subseteq I;$ then for any homogeneous elements $x,y\in B,$ 
$u\in I^\perp,$ and $v\in I,$ we get 
$b(xu,v)=-(-1)^{\bar{u}\bar{x}}b(u,xv)=0$ and $b([x,u,y],v)=
-(-1)^{\bar{x}\bar{u}}b([u,x,y],v)=(-1)^{\bar{x}(\bar{u}+\bar{y}+\bar{v})}b(u,[y,v,x])=0$ by (\ref{inv2}).
Then $BI^\perp\subseteq I^\perp$ and $[B,I^\perp,B]\subseteq I^\perp$ i.e.
 $I^\perp$ is an ideal of $B.$  \hfill $\square$\\
 
 We can now prove the following result.
 \begin{Theorem}\label{ThIperp}
 Let $B=B_0\oplus B_1$ be a Bol superalgebra such that $str(D_{x,y}L_z)=0$ 
 for all $x,y,z\in B.$ Then the Killing-Ricci form $\beta$ of $B$ is 
 nondegenerate if and only if the standard enveloping superalgebra 
 $L(B)=B\oplus IPS(B,B)$ is a semisimple Lie superalgebra.
 \end{Theorem}
{\bf Proof.} Let $\alpha$ be the Killing form of the Lie superalgebra 
$L(B).$ Then for all $x,y,z\in B,$ $\alpha(D_{x,y},z)=str(L_{D_{x,y}}L_z)=
str(D_{x,y}L_z)$ that is $str(D_{x,y}L_z)=0$ if and only if 
\begin{eqnarray}
\alpha(D_{x,y},z)=0 \label{cnes}
\end{eqnarray}
Then by (\ref{cnes}) and the invariance of $\alpha,$ we have for 
 $x,y,u,v\in B_0\cup B_1,$ 
 \begin{eqnarray}
 \alpha(D_{x,y},D_{u,v})&=&\alpha([x,y],[u,v])\nonumber\\
                        &=&-(-1)^{\bar{x}\bar{y}}\alpha(y,[x,[u,v]])\nonumber\\
                        &=&(-1)^{\bar{x}(\bar{u}+\bar{v}+\bar{y})}\alpha(y,    
                               [[u,v],x])\nonumber\\
                        &=&(-1)^{\bar{x}(\bar{u}+\bar{v}+\bar{y})}\beta(y,    
                               [u,v,x])\label{condNondeg}                    
 \end{eqnarray}
Thus, if $\beta$ is nondegenerate, the restriction of $\alpha$ to $IPS(B,B)\times IPS(B,B)$ is nondegenerate and therefore $\alpha$ is nondegenerate.
 
 Now by contradiction, suppose that $\beta$ is degenerate. Then by 
 Lemma \ref{Iperp}, $B^\perp$ is a nonzeo ideal of $B$ and therefore $B^\perp\oplus IPS(B,B^\perp)$ is a nonzero ideal of $L(B)$ by Lemma \ref{Idnvelop}.
 
 By the identities (\ref{cnes}) and (\ref{condNondeg}), we obtain
 \begin{eqnarray}
&&\alpha(B^\perp\oplus IPS(B,B^\perp),B\oplus IPS(B,B))\nonumber\\
 &=&\alpha(B^\perp,B)+\alpha(B^\perp, IPS(B,B))+\alpha(IPS(B,B^\perp),B)+
 \alpha(IPS(B,B^\perp),IPS(B,B))\nonumber\\
 &=&\alpha(B^\perp, B)+\alpha(B,[B^\perp, B, B])=0\nonumber
 \end{eqnarray}
It follows that $\alpha$ is nondegenerate  and therefore $L(B)$ 
 is nondegenerate, which ends the theorem. \hfill $\square$

\end{document}